\documentclass{article}

\usepackage{graphicx}

\usepackage{amssymb, verbatim, amsmath,amsfonts,natbib}
\usepackage{color}

\newtheorem{theorem}{Theorem}

\newtheorem{proof}{Proof}

\newcommand{\ds}{\displaystyle}

\newcommand{\tr}{\text{tr}}

\newcommand{\Perp}{\perp \! \! \! \perp}
\begin{document}

\title{Identifying Graphical Models
\thanks{The work was supported by a grant from the Swiss National Science
  Foundatation.}}

\author{Maya Shevlyakova and Stephan Morgenthaler \\
\'Ecole Polytechnique F\'ed\'erale de Lausanne\\
EPFL -- MATH-AA\\ Station 8\\ 1015 Lausanne, Switzerland\\\\
stephan.morgenthaler@epfl.ch 
}

\date{September 2012}
\maketitle

\begin{abstract}
  The ability to identify reliably a positive or negative partial
  correlation between the expression levels of two genes is influenced by
  the number $p$ of genes, the number $n$ of analyzed samples, and the
  statistical properties of the measurements.  Classical statistical theory
  teaches that the product of the root sample size multiplied by the size
  of the partial correlation is the crucial quantity. But this has to be
  combined with some adjustment for multiplicity depending on $p$, which
  makes the classical analysis somewhat arbitrary.  We investigate this
  problem through the lens of the Kullback-Leibler divergence, which is a
  measure of the average information for detecting an effect.  We conclude
  that commonly sized studies in genetical epidemiology are not able to
  reliably detect moderately strong links.  

\end{abstract}

\section{Introduction}
\label{intro}

Probabilistic graphical models are graphs in which nodes represent random
variables $X_u$ and the edges represent conditional dependence. Any two
nodes or variables that are not connected are independent, conditional on
the values of all the other random variables. Such models provide a compact
representation of a joint probability distribution. In a typical genetical
epidemiology application, the variables $X_u$ are gene expressions and
their measurements are available from tissue samples of $n$ patients. The
graphical model is used to describe the association between genes.  We
write $X_1 \Perp X_2| X_3 \ldots X_p$ to indicate that $X_1$ and $X_2$ are
conditionally independent, given $X_3 \ldots X_p$.  For the multivariate
normal distribution, conditional independence is equivalent to zero entries
in the inverse covariance matrix $\Sigma^{-1}$ (also called a concentration
or precision matrix). Thus, if $X$ is a $p$-dimensional normal random
vector with regular covariance matrix $\Sigma$, then for $1 \le u, v \le p$
with $u \neq v$
$$ X_u \Perp X_v \hspace{5pt}|\hspace{5pt} X_{rest}
\Longleftrightarrow \sigma^{uv} = 0\,,$$
where $\Sigma^{-1}= (\sigma^{uv})_{u,v = 1}^p$.

Estimating the structure of the concentration matrix from data can be
solved with a variety of statistical procedures. A possible approach for
low-dimensional data, for example,
consists in testing the inclusion of every edge separately, edge by
edge. Thus, we have to test $H_0^{uv}: \rho_{uv \cdot \text{rest}} = 0$ for
all $\binom{p}{2}$ choices of $u$ and $v$, where rest refers to the
variables with indices in $\{1,2,\ldots,p\}\setminus\{u,v\}$. To study the
feasibility of identifying the correct model, we could then
investigate how the power of this multiple testing problem depends on $n$
and $p$.

A better way to investigate the feasibility of edge-detection is based on
the Kullback-Leibler divergence. The Kullback-Leibler divergence (KLD)
measures in a statistically meaningful way the difference between two
probability distributions $F_1$ and $F_2$ with densities $f_1$ and $f_2$
\citep[see][]{kull}. It is defined as
\begin{equation}\label{KL_div}
D(f_1 | f_0) = \mathbb{E}_{X\sim F_1}\left[
\log\left(\frac{f_1(X)}{f_0(X)}\right)\right]\,.
\end{equation}
The divergence is thus simply the expected value of the
log-likelihood-ratio for a single observation $X$ from the alternative
model $F_1$ when testing the null model $F_0$. It is easy to show that this
divergence is positive unless $F_0=F_1$, in which case it is
zero. Furthermore, the bigger the KLD, the easier it is to distinguish
$F_1$ from $F_0$ by likelihood tests and the more powerful the test will be
\citep[for details, see][]{MorgStaudte2012}. If we dispose of $n$ independent
observations, the KLD is multiplied by $n$. If we test the absence of
partial correlations vs. the presence of partial correlations and assume
multivariate Gaussianity, the KLD is a useful tool to determine the average
amount of information in the data. Because it is based on likelihoods
rather than estimates, the KLD can be computed for any two models, without
reference to additional conditions such as $n > p$. This is an advantage of
this approach.

In the remainder of the paper, we will examine how information accumulates
when trying to fit a graphical model. When testing for edges, we will be
interested in the power of the test and the traditional asymptotic analysis
is not valid when $p > n$, while in the KLD approach, we can directly
compute the relevant amount of information.

\section{The Kullback-Leiber divergence}
\label{sec:1}

Suppose we have two $p$-variate normal populations with densities
\begin{equation}\label{normal}
f_i(x_1,\ldots,x_p) = \frac{1}{\ds|2\pi \Sigma_i|^{1/2}}
\exp (-\tfrac{1}{2}(x-\mu_i)^T\Sigma_i^{-1}(x-\mu_i))\text{ for }i=0,1\,,
\end{equation}
where $x=(x_1,\ldots,x_p)$ and $\mu_i,\Sigma_i$ denote the multivariate
means and covariance matrices. It follows that
\begin{multline*}\log \left(\ds\frac{f_1(x)}{f_0(x)}\right) =
\tfrac{1}{2}\log\left(\frac{|\Sigma_1|}{|\Sigma_0|}\right) -
\tfrac{1}{2}\tr\Sigma_0^{-1}(x-\mu_0)(x-\mu_0)^T \\
+\tfrac{1}{2}\tr\Sigma_1^{-1}(x-\mu_1)(x-\mu_1)^T\,.
\end{multline*}
Taking the expectation of the above, we can evaluate~(\ref{KL_div}) as
\begin{equation}\label{KL_normal}
D(f_1 | f_0) = \tfrac{1}{2}\log\left(\frac{|\Sigma_1|}{|\Sigma_0|}\right)+
\tfrac{1}{2}\tr\Sigma_0(\Sigma_1^{-1} -\Sigma_0^{-1})+
\tfrac{1}{2}\tr\Sigma_1^{-1}(\mu_0 - \mu_1)(\mu_0 - \mu_1)^T\,.
\end{equation}
We will make use of this formula for our purpose in which normal
populations with equal means but unequal covariance matrices are
compared. The null model will have a covariance matrix equal to the
identity matrix, while the alternative model will have a covariance matrix
whose inverse is nearly equal to the identity matrix. This describes a
situation where the $p$ variables have equal variance and only a very small
proportions of all partial correlations are non-zero.

\subsection{Divergence for a single non-zero partial correlation with known
placement}
\label{subsec:1}

Let $f_{uv}$ and $f_0$ be $p$ dimensional multivariate Gaussian densities
with mean $0$ and variances $\Sigma_{uv}$, $I_p$ (identity matrix ),
respectively. For now, the matrix $\Sigma_{uv}^{-1}$ has diagonal elements
equal to 1 and all off-diagonal values are zero, except for a value of
$\rho\neq 0$ in positions $(u,v)$ and $(v,u)$, that is, exactly one partial
correlation is non-zero. It is easy to show that the partial correlation is
equal to $-\rho$ in this case. We will write $\Sigma_{uv}^{-1} = I_p +
U_{uv}$, where $U_{uv}$ contains the off-diagonal elements of
$\Sigma_{uv}$.  To compute the divergence (\ref{KL_normal}), we need the
determinant of $\Sigma_{uv}^{-1}$ and the trace of $\Sigma_{uv}$, which are
$$|\Sigma_{uv}^{-1}| = 1 + (-1)^{u+v} \rho (-1)^{v-1 + u} \rho
= 1+\rho^2(-1)^{2u+2v-1} = 1-\rho^2$$ and $\text{trace}(\Sigma_{uv})=p-2 +
2/(1-\rho^2)$. Note that the determinant does not depend on the
position $(u,v)$, nor does it depend on the dimension $p$, while for the
trace the diagonal elements of $\Sigma_{uv}$ are needed. They are are equal
to 1, except in positions $(u,u)$ and $(v,v)$, where they are
$1/(1-\rho^2)$. From the previous expression (\ref{KL_normal}) we then find
\begin{equation}
D(f_{uv}|f_0) =
\frac{\log(1-\rho^2)}{2}+ \frac{\rho^2}{1-\rho^2}\,.\label{divergence}
\end{equation}
A sample of $n$ observations drawn from the alternative model $f_{uv}$ has
an information content in favor of rejecting the null model $f_0$ of 
$nD(f_{uv}|f_0)$, which tends to infinity as $n$ grows larger. This is true
for any value $\rho\neq 0$. For a large enough sample, even a slight partial
correlation between the u-th and v-th variable will for sure be detected.

\subsection{Divergence when the placement is unknown}\label{subsec:2}

In our formulation of the density $f_{uv}$, we of course make use of the
knowledge of the placement of the positive partial correlation. Because of
this, Eq.~(\ref{divergence}) is only useful in understanding the null
hypothesis $H_0^{uv}$, which is only one of the possibilities one has to
examine in practice. How does the divergence change if we do not know the
pair of correlated variables? To answer this question, we consider a
different alternative model, namely the multivariate mixture density
\begin{equation}\label{Normal_mixture}
f_1(x_1,\ldots,x_p) = 
\sum_{u=1}^{p-1}\sum_{v=u+1}^p \frac{1}{\binom{p}{2}}f_{uv}(x_1,\ldots,x_p)\,.
\end{equation}
In this case, the following result holds.
\begin{theorem} The amount of information about $\rho$ in a sample of size
  $n$ drawn from the mixture distribution (\ref{Normal_mixture}),
  information in favor of distinguishing this mixture model from the null
  model of independent standard normals, is for $\vert\rho\vert <
  1/\sqrt{2}$ and large dimension $p$ equal to
\begin{equation}
 nD(f_1 | f_0)=\frac{4n}{p}\left(\frac{1-\rho^2}{\sqrt{1-2\rho^2}}-1\right)+
   o(p^{-1})\,.\label{eq:approxKLD}
\end{equation}

\end{theorem}
\begin{proof}
Since 
$$f_{uv}(x_1,\ldots,x_p) = \frac{\sqrt{1-\rho^2}}{(2\pi)^{p/2}}
\exp\left(-\frac{1}{2}(x_1^2 + \ldots +x_p^2+ 2\rho x_ux_v)\right)\,,$$
the likelihood ratio is equal to
\begin{multline*}
\frac{f_1(x)}{f_0(x)} =
\frac{\ds\frac{\sqrt{1-\rho^2}}{{p\choose 2}} \sum_{u=1}^{p-1}\sum_{v=u+1}^p
\exp\left(-\tfrac{1}{2}(x_1^2 + \ldots +x_p^2+ 2\rho x_ux_v)\right)}
{\ds\exp(-\tfrac{1}{2}(x_1^2 +\cdots+x_p^2))}\\
={\ds\frac{\sqrt{1-\rho^2}}{{p\choose 2}}
\sum_{u=1}^{p-1}\sum_{v=u+1}^p\exp (-\rho {x_u}{x_v})}\,,
\end{multline*}
where as before $x=(x_1,\ldots,x_p)$.
It follows that
\begin{multline}\label{KL_div_mix}
D(f_1 | f_0) = \int_{\mathbb{R}^p} \log\left(\frac{f_1(x)}{f_0(x)}\right)f_1(x)\,dx \\=
\int_{\mathbb{R}^p} \log\left({\ds\frac{\sqrt{1-\rho^2}}{{p\choose 2}}
\sum_{u=1}^{p-1}\sum_{v=u+1}^p\exp (-\rho{x_u}{x_v})}\right)f_1(x) \, dx\\=
\int_{\mathbb{R}^p} \log\left({\ds\frac{\sqrt{1-\rho^2}}{{p\choose 2}}
\sum_{u=1}^{p-1}\sum_{v=u+1}^p\exp (-\rho{x_u}{x_v})}\right)f_{12}(x) \, dx\,,
\end{multline}
where the last equality follows from the fact that the ratio of the densities
is invariant with respect to permutations of the components of $x$.

The expectation (\ref{KL_div_mix}) can be approximated for large dimensions
$p$ via the law of large numbers. Let $(X_1,\ldots,X_p)$ be a random vector
with density $f_{12}$. It follows that if $(Z_1,\ldots,Z_p)$ are
independent unit Gaussian random variables, we have the representation
$(X_1=AZ_1+BZ_2,X_2=AZ_2 +BZ_1, X_3=Z_3,\ldots,X_p=Z_p)$, where
$A=(1+[(1+\rho)/(1-\rho)]^{1/2})/(2(1+\rho)^{1/2})$ and 
$B=(1-[(1+\rho)/(1-\rho)]^{1/2})/(2(1+\rho)^{1/2})$. The
integrand in (\ref{KL_div_mix}) can thus be written as
\begin{multline}
\log\left( \frac{\sqrt{1-\rho^2}}{{p\choose 2}}\left(
e^{-\rho{X_1}{X_2}}+ \sum_{v=3}^p\left(e^{-\rho X_1Z_v}+e^{-\rho X_2Z_v}\right)+
\sum_{u=3}^{p-1}\sum_{v=i+1}^pe^{-\rho Z_uZ_v}\right)\right)
\\=
\log\Bigg(\sqrt{1-\rho^2} \Bigg(O_P(p^{-2})+\frac{2}{p^2-p}
\left(\sum_{v=3}^p\left(e^{-\rho X_1Z_v}+e^{-\rho X_2Z_v}\right)\right)
\\+
\frac{2}{p^2-p}\left(\sum_{u=3}^{p-1}\sum_{v=i+1}^pe^{-\rho Z_uZ_v}
\right)\Bigg)\Bigg)\,.\label{eq:KLD}
\end{multline}
For large $p$, the two last terms can be approximated by their asymptotic
limits. The product of two independent normal variables, which appears in these
expressions, has surprising properties \citep[see][]{aroian1947}. Elementary
calculations show that if $Z_1,Z_2$ are independent with a unit
Gaussian distribution, then $\mathbb{E}[\exp(-\rho Z_1Z_2)]=(1-\rho^2)^{-1/2}$.
With $Y=\exp(-\rho Z_1Z_2)$ we thus find for $\gamma > 0$ that
$\mathbb{E}[\vert Y\vert^\gamma] = \mathbb{E}[\exp(-\gamma\rho Z_1Z_2)]=
(1-\rho^2\gamma^2)^{-1/2}$, which is finite if $\vert\rho\vert\gamma <
1$. For $0 < \vert\rho\vert < 0.5$, this allows a maximal value $\gamma
\geq 2$ (finite variance), while for $0.5 \leq \vert\rho\vert < 1$, we
always have a maximal value $1<\gamma < 2$ (infinite variance). 

With regard to the last term in (\ref{eq:KLD}), this shows that for all
$\rho\neq 0$ we have convergence in probability of the mean
\begin{multline*}
\hspace{1cm}\binom{p}{2}^{-1}\left(\sum_{u=3}^{p-1}\sum_{v=i+1}^pe^{-\rho Z_uZ_v}
\right)=\\\left(1-\frac{4}{p}+O(p^{-2})\right)
\binom{p-2}{2}^{-1}\left(\sum_{u=3}^{p-1}\sum_{v=i+1}^pe^{-\rho Z_uZ_v}
\right)\overset{\text{in P}}{\longrightarrow} \frac{1}{\sqrt{1-\rho^2}}
\text{ as }p\to\infty\,,
\end{multline*}
but the rate of convergence in the law of large numbers depends on the value
of $\rho$. 

The other term in (\ref{eq:KLD}) involves either $X_1$ or $X_2$, which both
have variance $(1-\rho^2)^{-1}$ and have the same marginal distribution as
$Z_1(1-\rho^2)^{-1/2}$. It follows that the typical summand is $Y =
\exp(-\rho(1-\rho^2)^{-1/2} Z_1Z_2)$, which for $\gamma > 0$ satisfies
$\mathbb{E}[\vert Y\vert^\gamma] =
\mathbb{E}[\exp(-\gamma\rho(1-\rho^2)^{-1/2} Z_1Z_2)]=
(1-\rho^2\gamma^2/(1-\rho^2))^{-1/2}$. This is only finite, if $\gamma^2 <
(1-\rho^2)/\rho^2$. For $\vert \rho\vert < 0.5$, $\gamma$ takes values up
to at least $\sqrt{3}$ and for $\vert \rho\vert < 1/\sqrt{5}$, the variance
is finite, that is, $\gamma$ is at least 2. For $\vert\rho\vert > 1/\sqrt{2}$,
$\gamma$ is less than 1 and the expected value becomes infinite. In the
range of $\rho$ values we are considering, the expectation is equal to 
$(1-\rho^2/(1-\rho^2))^{-1/2}=[(1-\rho^2)/(1-2\rho^2)]^{1/2}$ and thus
\begin{eqnarray*}
\binom{p}{2}^{-1}\left(\sum_{v=3}^{p}e^{-\rho X_1Z_v}
\right)=\left(\frac{2}{p}+O(p^{-2})\right)
\frac{1}{p-2}\left(\sum_{v=3}^{p}e^{-\rho X_1Z_v}
\right)
\end{eqnarray*}
and
\begin{eqnarray*}
\frac{1}{p-2}\left(\sum_{u=3}^{p}e^{-\rho X_1Z_v}
\right)
\overset{\text{in P}}{\longrightarrow} \sqrt{\frac{1-\rho^2}{1-2\rho^2}}\text{ as }p\to\infty\,.
\end{eqnarray*}
Substitution of the limit as an approximate value for large $p$ leads to
the following expansion of the value inside the logarithm of (\ref{eq:KLD})
\begin{equation*}
\sqrt{1-\rho^2}\left(\frac{4}{p}\frac{\sqrt{1-\rho^2}}{\sqrt{1-2\rho^2}}+ \left(1-\frac{4}{p}
  \right)\frac{1}{\sqrt{1-\rho^2}}\right)\\
=1-\frac{4}{p} +\frac{4}{p}\frac{1-\rho^2}{\sqrt{1-2\rho^2}} \,.
\end{equation*}
Expanding the logarithm to the required order leads to the approximate KLD
value claimed in the theorem.
\end{proof}

\section{Discussion and extensions}

\subsection{Other values of $\rho$ and numerical comparisons} 

As the value of $\vert\rho\vert$ increases towards 1, the approximate
computation of the KDL undergoes several transitions. For the main term in
(\ref{eq:KLD}) the only point of transition occurs at $\rho = 0.5$, when
the variance becomes infinite. For the minor term they occur at $\rho =
1/\sqrt{5}$, when the variance becomes infinite and at $\rho=1/\sqrt{2}$,
when the expectation becomes infinite. At this second point, our formula is
no longer valid, because the terms of order $p^{-1}$ are not the leading
terms. In this case, a more refined analysis of the tail probabilities of
the law of $Y$ is required \citep[see, for example,][]{gut2004,baum1965,
  feller1945}.

Figure 1 includes the numerical results for $\rho=0.9$. The plot makes it
clear that in this case the information content decreases less rapidly with
increasing dimension $p$. Even in this case, the linearity in the plot of
log(KLD) as a function of log(p) remains, but the slope passes from $-1$ to
$-0.25$. The analysis based on moments of $Y$ suggests for $\vert\rho \vert
> 1/\sqrt{2}$ an order of $p^{-\tau}$ with $0< \tau < 2(1-\rho^2)$
\citep[see][Theorem 1]{baum1965}.

\begin{figure}
\begin{center}
\includegraphics[scale = 0.6]{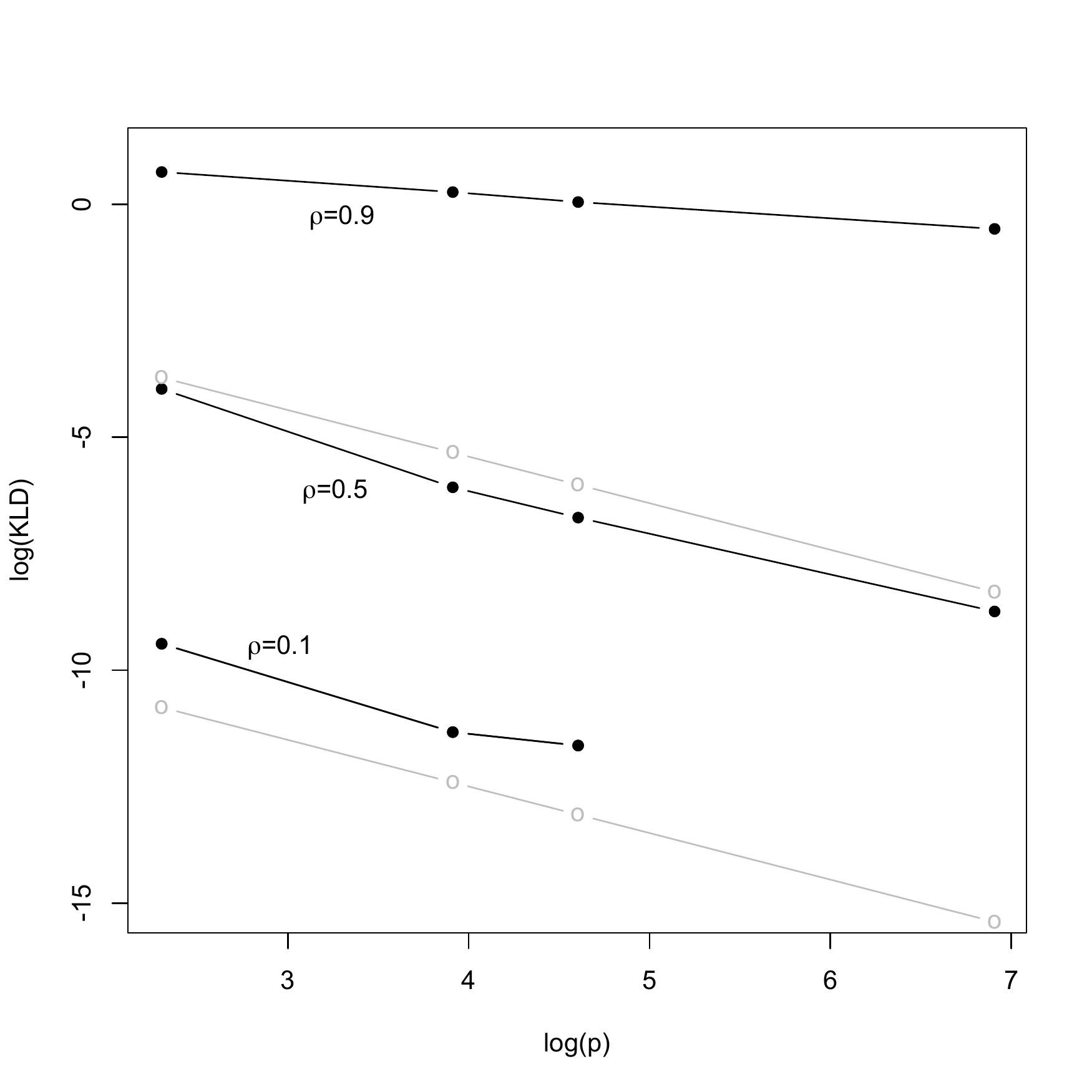}
\caption{The comparison of the numerical (black) and the approximate (grey)
  divergence values show a good agreement, even for $p=10$. The approximate
  divergence decreases as 1/p, which translates into a slope of $-1$ in our
  plot.}
\label{fig:1}      
\end{center}
\end{figure}

Figure~\ref{fig:1} compares the values obtained by Monte Carlo
sampling with the approximation given in Theorem 1. The agreement is quite
good.  

\subsection{Classical asymptotics} 

The analysis using the KLD is related, but different, from the more widely
known asymptotic or local power. When using the KLD, there is no correction
for multiplicity involved, no constraints of the type $n>p$ are needed and
no limits towards infinite study sizes are taken. The KLD thus provides a
more solid foundation for judgeing the sample sizes needed in order to
reliably detect effects. Here we briefly compare it with the traditional
asymptotic approach, where $n\to\infty$ (and implicitly $n > p$).  When
testing the null hypothesis $H_0^{uv}: \rho_{uv \cdot \text{rest}} = 0$
against one-sided alternatives $\rho_{uv \cdot \text{rest}} > 0$ based on
the estimator $s_{uv \cdot \text{rest}}$ and the Bonferroni correction for
the number of tests $m=p(p-1)/2$, the power function for large sample sizes
is approximately equal to
\begin{equation}
1-\Phi\left(z_{1-\alpha/m}-\sqrt{n}\rho_{uv \cdot \text{rest}}
  \right)=\Phi\left(\sqrt{n}\rho_{uv \cdot \text{rest}}
  -z_{1-\alpha/m}\right)\,,\label{eq:asypow}
\end{equation}
where $z_{1-\alpha/m}$ denotes the $1-\alpha/m$ quantile of the standard
normal distribution. Using the asymptotic approximation for this quantile
$z_{1-\alpha/m}\sim \sqrt{2\log(m/\alpha)}$, leads to the following
one-sided local asymptotic power at the alternative $\rho_{uv \cdot \text{rest}}>0$:
\begin{equation}
\Phi\left(\sqrt{n}\rho_{uv \cdot \text{rest}}-
\sqrt{2\log\left[p(p-1)/(2\alpha)\right]}\right)\,,\label{eq:asyp}
\end{equation} 
which depends on $p$ via the logarithm.
The above approximation of the quantile $z_{1-\alpha/m}$
is quite crude and gives values that are typically too large so that the
power might be underestimated.  This local asymptotic approximation is
based on the asymptotic normality of the estimator of the partial
correlation and on the consideration of alternatives close to the null
hypothesis \citep[see for example Chapter 10 of][]{serfling}. Finding the
approximate power involves the calculation of the slope and is based on the
expected value of the partial correlation estimate. This can be shown to be
\citep[equation (18), section 5.1 in][]{muirhead1982aspects}
$$\frac{2}{f}\,(\Gamma[(f+1)/2]/\Gamma[f/2])^2\,\rho_{uv \cdot \text{rest}}
~_2F_1\left[\tfrac{1}{2},\tfrac{1}{2};(f+2)/2,\rho^2_{uv \cdot
    \text{rest}}\right],$$ 
where $f=n+1-p$ and $~_2F_1[\cdot]$ is a hypergeometric function, which is
evaluated at $\rho_{uv \cdot \text{rest}}$. It follows that its
derivative with respect to $\rho_{uv \cdot \text{rest}}$, evaluated at
$\rho_{uv \cdot \text{rest}}=0$, is equal to 
$2(\Gamma[f/2])^2/(f\Gamma[(f+1)/2]^2)$, which tends to 1 as $f\to\infty$
by Stirling's approximation. Because
the asymptotic variance of the partial correlation estimator assuming that
the null hypothesis is true is $1$, the slope of the test or its Pitman
efficacy is equal to 1.  

A comparison between (\ref{eq:approxKLD}) and (\ref{eq:asyp}) can be
based on the fact that in order to reach a power of about 0.5 at level
$\alpha$, the KLD of an experiment must exceed $z_{1-\alpha}^2$
\citep[see][]{MorgStaudte2012}. From this, one can derive a formula for the
needed size of a study, $n =
pz_{1-\alpha}^2/[4((1-\rho^2)/\sqrt{1-2\rho^2}-1)]$. The equivalent value
of $n$ from the asymptotic power on the other hand predicts that
$n=(z_{1-\alpha/m}/\rho)^2$, where $m=p(p-1)/2$ is the number of tests. For
values of $\rho<0.5$, the KLD-based formula gives much higher values of the
study size $n$. For example, around $n=20,000$ subjects would be required
to detect a partial correlation in a single pair of $p=1000$ genes. The
asymptotic power wrongly suggests that $n=135$ subjects would be
sufficient.  Generally speaking, when $\rho<0.5$ the problem of identifying
a partial correlation is hopeless, unless the number of candidate genes
that are tested can be reduced below $p=100$. Figure 1 also gives an
indication of what will happen for a strong effect, $\rho=0.9$. The value
of KLD decreases by about a factor of 0.24 for each increase of $p$ by a
factor of 10. If we extrapolate to $p=10^6$, we have a KLD value of about
$8\times10^{-3}$. We thus would need a study involving at least $n=330$
subjects, which is doable.

\subsection{Detecting correlations}

The model in which the covariance matrix is equal to our precision matrix,
that is, the model with $p$ measurements with equal variance and a single
non-null covariance $\rho$ has been analyzed by \cite{arias2012}. In
computations not shown here, we obtain the following result which holds for
small values of $\rho$ and large values of $p$
\begin{equation*}
 nD(f_1 | f_0)=\frac{2n}{p^2}\left(\frac{1}{\sqrt{1-\rho^2}}-1\right)+
   o(p^{-2})\,.
\end{equation*}
Since the order of the leading term is $p^{-2}$, the detection is this
model is nearly impossible unless $n$ is very large.

\subsection{Divergence for two partial correlations }

\begin{figure}
\begin{center}
\includegraphics[scale = 0.6]{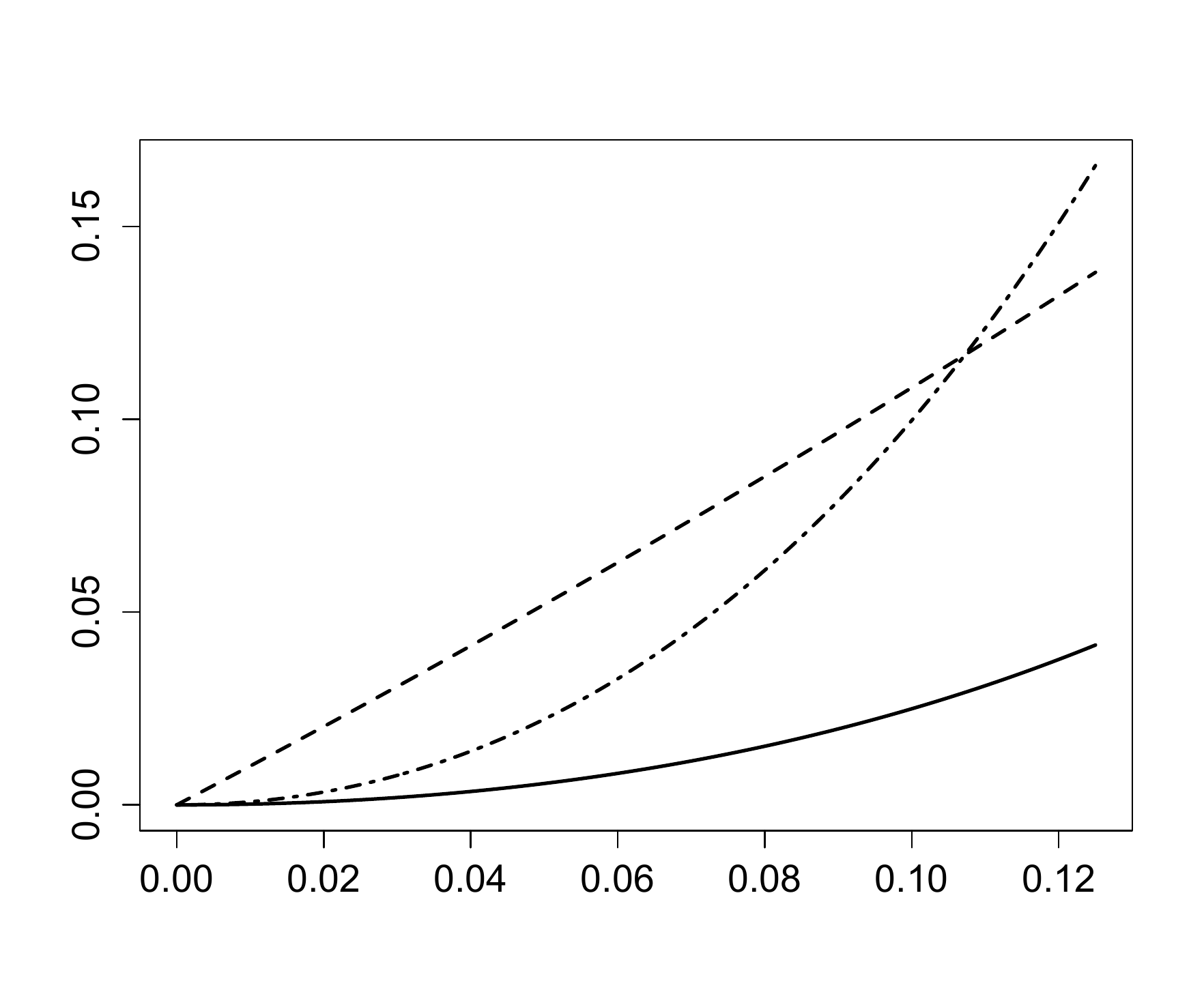}
\caption{The qualitative behavior of the Kullback-Leibler divergence in the
three cases discussed in this paper are shown here. The horizontal axis is
equal to the square $\rho^2$. The solid curve is the
case of a change a single off-diagonal element in the precision matrix,
while the dashed curve is a single off-diagonal element in the correlation
matrix. In both cases, the diagonal elements of the matrix are all equal to
one. The dotted/dashed curve is the case of two non-null and equal elements
in the precision matrix.\label{fig:three}}
     
\end{center}
\end{figure}

If the we consider the alternative multivariate Gaussian model with an
inverse covariance matrix in which the diagonal elements are 1 and exactly two
pairs of variables have a partial correlation of $\rho$, then the following
result holds for  $\vert\rho\vert < 1/\sqrt{2}$
\begin{equation*}
nD(f_1 | f_0)=\frac{16n}{p}\left(\frac{1-\rho^2}{\sqrt{1-2\rho^2}}-1\right)
+o(p^{-1})\,.
\end{equation*}
Figure \ref{fig:three} shows the qualitative behavior of the three
functions we computed. Note that in the case of a perturbation of the
correlation matrix by a single non-null element, the information increases
very slowly and is of order $O(n/p^2)$, while a perturbation on the level
of the precision matrix leads to $O(n/p)$. If two couples are correlated
with an equal correlation rather than a single couple, the information gain
is four-fold.

\bibliographystyle{chicago}
\bibliography{kld}   

\end{document}